\documentclass[11pt,a4paper]{amsart}

\usepackage{graphicx}
\usepackage{amssymb}

\usepackage[unicode,hidelinks,bookmarks]{hyperref}
\hypersetup{
colorlinks=true,
citecolor=[rgb]{0,0,0.5},
linkcolor=[rgb]{0,0,0.5},
urlcolor=[rgb]{0,0,0.75},
pdfpagemode=UseNone,
pdfstartview=FitH,
pdfdisplaydoctitle=true,
pdftitle={Band-limited maximizers on the circle},
pdfauthor={Oliveira e Silva, Thiele, Zorin-Kranich},
pdflang=en-US
}

\def\R{\mathbb{R}}

\def\Z{\mathbb{Z}}

\newcommand{\pvector}[1]{
  \begin{pmatrix}
    #1
  \end{pmatrix}} %
\newcommand{\ddirac}[1]{
  \,\boldsymbol{\delta}\!\pvector{#1}\!} %

\renewcommand{\d}{\,{\rm d}}
\renewcommand{\k}{\mathbf{k}}

\newtheorem{theorem}{Theorem}

\newtheorem{lemma}[theorem]{Lemma}

\title[Band-limited maximizers on the circle]{Band-limited maximizers for a Fourier extension inequality on the circle}

\author[Oliveira e Silva]{Diogo Oliveira e Silva}
\author[Thiele]{Christoph Thiele}
\author[Zorin-Kranich]{Pavel Zorin-Kranich}

\address{
School of Mathematics\\
University of Birmingham\\ 
Birmingham, B15 2TT, England}
\email{d.oliveiraesilva@bham.ac.uk}

\address{
        Hausdorff Center for Mathematics\\
        53115 Bonn, Germany}
\email{thiele@math.uni-bonn.de}
\email{pzorin@math.uni-bonn.de}

\begin{document}
\begin{abstract}
Among the class of functions with Fourier modes up to degree 30, constant functions are the unique real-valued maximizers for the endpoint Tomas--Stein inequality on the circle.  
\end{abstract}

\maketitle

\section{Introduction}

We are interested in the sharp constant of the endpoint Tomas--Stein adjoint restriction inequality \cite{To75} on the circle $\mathbb{S}^1=\{\omega\in\R^2:|\omega|=1\}$. More precisely, we seek 
a maximizer for the functional $\Phi$ defined on nonzero $f\in L^2(\mathbb{S}^1)$ by 
$$\Phi(f) := \|\widehat{f\sigma}\|^6_{L^6(\R^2)}\|f\|_{L^2(\mathbb{S}^1)}^{-6}.$$
We have written $\sigma$ for the arc length measure on the circle $\mathbb{S}^1$, and we have used the Fourier transform  
\begin{equation*}
\widehat{f\sigma}(x):=\int_{\mathbb{S}^1} f(\omega)\,e^{-i x\cdot\omega}\, \d\sigma_\omega,\;\;\;(x\in\R^2).
\end{equation*}
It is known that maximizers of $\Phi$ exist \cite{Sh16} and are smooth \cite{Sh16b}, and that the constant function ${\bf 1}$ is a local maximizer of $\Phi$, see \cite[Theorem 1.1]{CFOST17}. 
Moreover, real-valued maximizers of $\Phi$ are known to be nonnegative, antipodally symmetric functions, that is 
$$f(\omega)\ge 0,\ f(\omega)=f(-\omega),$$ 
for every $\omega\in \mathbb{S}^1$. 
It is natural to conjecture that constant functions are global maximizers of $\Phi$, in which case a complete characterization of complex-valued maximizers is given by \cite[\S1, Step 6]{CFOST17}.

In this paper, we report on numerical verification of a finite dimensional variant of this conjecture:

\begin{theorem}\label{thm:main}
Let $f\in L^2(\mathbb{S}^1)$ be non-negative and antipodally symmetric. 
Assume that $\widehat{f}(n)=0$ if $|n|>30$. Then
$$\Phi(f)\leq\Phi({\bf 1}),$$
with equality if and only if $f$ is constant.
\end{theorem}

A numerical difficulty in this problem is that there are close competitors for maximizers, namely functions that concentrate in the vicinity of two antipodal points. 
Heisenberg uncertainty allows for functions with Fourier modes up to degree $30$ to localize roughly in a $\frac{2\pi}{30}$-neighborhood of these points.

This paper continues efforts to implement Foschi's program \cite{Fo15}
for the 2-sphere in the case of the circle, see also \cite{CFOST17}. The
approach works through positive semi-definiteness of a certain quadratic form on the relevant finite dimensional space. It would be nice to establish positive semi-definiteness for the full space. For recent similar work on the paraboloid, \mbox{see \cite{Go17}.}

\section{Proof of Theorem \ref{thm:main}}
 
 With  $f$ as in  Theorem \ref{thm:main}, we compute  
	$$\|\widehat{f\sigma}\|_{L^6(\R^2)}^6
	= (2\pi)^2 \int_{(\mathbb{S}^1)^6} 	\ddirac{\sum_{j=1}^6 \omega_j }   \Big(\prod_{j=1}^3 f_j(\omega_j)\d\sigma_{\omega_j}\Big)\Big(\prod_{j=4}^6 f_j(-\omega_j)\d\sigma_{\omega_j}\Big) $$
 	$$= 5 \pi^2 \int_{(\mathbb{S}^1)^6} 	\ddirac{\sum_{j=1}^6 \omega_j }  (|\omega_4+\omega_5+\omega_6|^2-1) 
 	\Big(\prod_{j=1}^3 f_j(\omega_j)\d\sigma_{\omega_j}\Big)
 	\Big(\prod_{j=4}^6 f_j(-\omega_j)\d\sigma_{\omega_j}\Big)
    $$
  	$$\le 5\pi ^2 \int_{(\mathbb{S}^1)^6} 	\ddirac{\sum_{j=1}^6 \omega_j}  (|\omega_4+\omega_5+\omega_6|^2-1) \Big(\prod_{j=1}^3 f_j(\omega_j)^2\Big)
  	\prod_{j=1}^6 \d\sigma_{\omega_j} $$
  	$$\le  5 \pi^2 \frac{\|f\|_{L^2(\mathbb{S}^1)}^6}{\|{\bf 1}\|_{L^2(\mathbb{S}^1)}^6} 
	\int_{(\mathbb{S}^1)^6} 	\ddirac{\sum_{j=1}^6 \omega_j}  (|\omega_4+\omega_5+\omega_6|^2-1)
\prod_{j=1}^6 \d\sigma_{\omega_j}= \Phi({\bf 1}) \|f\|_{L^2(\mathbb{S}^1)}^6.$$
 
Here the first line is simply Plancherel's identity. The second line 
is Foschi's idea to improve the situation by artificially
inserting a weight, which after symmetrization over the
indices $j$ reverts to a constant, see the computation following \cite[Lemma 1.3]{CFOST17}.
The third line is the crucial inequality. We defer its proof for the moment.
The inequality in the fourth line is an application of the main result proved in  \cite[Theorem 1.2]{CFOST17}. Equality is attained if and only
if $f$ is constant. Identification of the constant in the fourth line
is then easy and was also observed in \cite{CFOST17}.

This proves Theorem \ref{thm:main}, safe for verification of the crucial inequality in the third line. Note that this inequality would follow from
$$|\omega_1+\omega_2+\omega_3|=|\omega_4+\omega_5+\omega_6|$$
and the inequality between the arithmetic mean and the geometric mean,
$$\prod_{j=1}^3 f_j(\omega_j)\prod_{j=4}^6 f_j(-\omega_j)
\le \frac 12\Big(\prod_{j=1}^3 f_j(\omega_j)^2
+\prod_{j=4}^6 f_j(-\omega_j)^2\Big),$$
if the weight were positive. Unfortunately, the weight is not positive.
One reason to believe that the inequality  still holds as stated is that the negative part of the weight is small, and via antipodal symmetry the values of the function on the negative part of the weight have a strong correlation with the values on the positive part. 
However, the support of the measure	 $$\ddirac{\sum_{j=1}^6 \omega_j}$$
is not preserved under antipodal symmetry, which makes it difficult to exploit this correlation. We resort to numerical verification
of the crucial inequality in the given finite dimensional space of functions.

Consider the index set
$$Z=\{k\in 2\mathbb{Z}, |k|\le 30\},$$
and expand the band-limited function $f$ into a Fourier series
$$f(\omega)=\sum_{k\in Z} a_k \omega^k,$$
where we identify $\mathbb{R}^2$ with the complex plane and correspondingly define products and powers of elements in $\mathbb{R}^2$.
Note that 
\begin{equation}
\label{realsym}a_{-k}=\overline{a_k}
\end{equation} 
for every $k\in Z$.
We write a constant multiple of the left-hand side of the crucial \mbox{inequality as}
$$ \sum_{{\bf k}\in \Z^6}
 L_{\bf k} \Big(\prod_{j=1}^3 a_{k_j}\Big)\Big(\prod_{j=4}^6 a_{-k_j}\Big),
$$
$$L_{{\bf k}}:= (2\pi)^{-5}\int_{(\mathbb{S}^1)^6} 	\ddirac{\sum_{j=1}^6 \omega_j}   (|\omega_4+\omega_5+\omega_6|^2-1) \Big(\prod_{j=1}^3 \omega_j^{k_j}\Big)\Big(\prod_{j=4}^6 \omega_{j}^{-k_j}\Big)
\prod_{j=1}^6 \d\sigma_{\omega_j},$$
and the same multiple of the right-hand side as

$$ \sum_{{\bf k}\in \Z^6}
R_{\bf k} \Big(\prod_{j=1}^3 a_{k_j}\Big)\Big(\prod_{j=4}^6 a_{-k_j}\Big),$$
$$R_{{\bf k}}:= (2\pi)^{-5}\int_{(\mathbb{S}^1)^6} 	\ddirac{\sum_{j=1}^6 \omega_j}   (|\omega_4+\omega_5+\omega_6|^2-1) \Big(\prod_{j=1}^3 \omega_j^{k_j-k_{j+3}}\Big)
\prod_{j=1}^6 \d\sigma_{\omega_j}.$$

Define for ${\bf m}\in Z^3$
$$s_{\bf m}:=a_{m_1}a_{m_2}a_{m_3},$$ 
and note that $(s_{\bf m})_{{\bf m}\in Z^3}$ is an element of $Sym(Z^3)$, the vector space of functions on $Z^3$ symmetric under permutation of the three indices.
At this point we do not require the symmetry \eqref{realsym},
instead we pass to a larger space allowing for a convenient orthogonal splitting later.

The crucial inequality then follows from positive semi-definiteness
of the quadratic form 
$$\sum_{{\bf m},{\bf n}\in Z^3} 
Q_{{\bf m},{\bf n}}s_{\bf m} \overline{s_{ \bf n}}:=
\sum_{{\bf m},{\bf n}\in Z^3} \frac 16 \sum_{\sigma\in S_3}
(R_{{\bf m},{\bf n_\sigma}}-L_{{\bf m},{\bf n_\sigma }})s_{\bf m} \overline{s_{ \bf n}}$$
on $Sym(Z^3)$,
where we write $S_3$ for the group of permutations of three elements and
$${\bf n}_\sigma=(n_{\sigma(1)},n_{\sigma(2)},n_{\sigma(3)}).$$
Note that the symmetrization over $S_3$ does not change the value of the quadratic form whenever the coefficients $s_{\bf n}$ are symmetric. It merely symmetrizes the coefficients of the quadratic form, and allows to reduce the dimension of the matrix by identifying equivalent
tuples. Letting $\tilde{X}$ be the space of tuples in $Z^3$
satisfying $m_1\le m_2\le m_3$, it therefore suffices to prove 
positive definiteness of the quadratic form
$$Q(s,s)= \sum_{{\bf m},{\bf n}\in \tilde{X}} Q_{{\bf m},{\bf n}} s_{\bf m} \overline{s_{\bf m}} .$$

Note that the matrix 
 $(Q_{{\bf m},{\bf n}})_{{\bf m},{\bf n}\in \tilde{X}}$ is Hermitian.
Moreover, for all ${\bf m}\in Z^3$ we have
$$R_{{\bf m},(0,0,0)}=L_{{\bf m},(0,0,0)}$$
and hence the Dirac delta vector $\delta_{(0,0,0)}$ corresponding to constant functions on the circle is in the kernel
of the matrix $(Q_{{\bf m},{\bf n}})_{{\bf m},{\bf n}\in \tilde{X}}$.
Therefore we replace $\tilde{X}$ by $X = \tilde{X} \setminus \{(0,0,0)\}$.

 A change of variables 
$$(\omega_1,\omega_2,\omega_3,\omega_4,\omega_5,\omega_6)\mapsto 
(\omega_1\cdot\omega,\omega_2\cdot\omega,\omega_3\cdot\omega,\omega_4\cdot\omega,\omega_5\cdot\omega,\omega_6\cdot\omega)$$
for some arbitrary $\omega$ of modulus one in the expressions for $R_{\bf k}$ and $L_{\bf k}$ shows that
 $$Q_{{\bf m},{\bf n}}=\omega^{d({\bf m})-d({\bf n})}Q_{{\bf m},{\bf n}},$$
where we have denoted
$$d({\bf m})=m_1+m_2+m_3.$$
We conclude
 $$Q_{{\bf m},{\bf n}}=0$$
whenever $d({\bf m})\neq d({\bf n})$. The matrix $(Q_{{\bf m},{\bf n}})_{{\bf m},{\bf n}\in X}$
therefore has the structure of a diagonal block matrix, with blocks enumerated by $D:=d({\bf m})$.
It suffices to prove positive semi-definiteness for each block $(Q_{{\bf m},{\bf n}})_{{\bf m},{\bf n}\in X_{D}}$ separately, where $X_{D} = \{{\bf m}\in X : d({\bf m}) = D\}$.
This will be done in the following section.

\section{Numerical computations}\label{sec:Numerics}
In order to verify that the matrix $(Q_{{\bf m},{\bf n}})_{{\bf m},{\bf n}\in X_{D}}$ is positive definite, we split it into a numerically computed approximation and an error term.
The smallest eigenvalue of the numerical approximation will be larger than the operator norm of the error term.

 \begin{table}[ht]
 \centering
 \begin{tabular}{|c | c || c | c || c |c || c | c |}
 \hline
 $D$&$\lambda_{\min}$&$D$&$\lambda_{\min}$&$D$&$\lambda_{\min}$& $D$ &$\lambda_{\min}$\\[0.5ex]
 \hline
  0&0.00035 & 24&0.00061 & 48&0.00121 & 72&0.00407\\
  2&0.00037 & 26&0.00064 & 50&0.00133 & 74&0.00501\\ 
  4&0.00038 & 28&0.00067 & 52&0.00144 & 76&0.00596\\
  6&0.00042 & 30&0.00069 & 54&0.00154 & 78&0.00668\\
  8&0.00045 & 32&0.00073 & 56&0.00171 & 80&0.00937\\
 10&0.00049 & 34&0.00077 & 58&0.00188 & 82&0.01258\\
 12&0.00052 & 36&0.00081 & 60&0.00203 & 84&0.01332\\
 14&0.00055 & 38&0.00086 & 62&0.00229 & 86&0.02997\\
 16&0.00057 & 40&0.00092 & 64&0.00255 & 88&0.04400\\
 18&0.00057 & 42&0.00097 & 66&0.00278 & 90&0.20081\\
 20&0.00058 & 44&0.00105 & 68&0.00324  & &\\
 22&0.00059 & 46&0.00113 & 70&0.00369 & &\\[1ex]                          
 \hline
 \end{tabular}
 \caption{Minimal eigenvalue for the approximation for the block 
 	$D\in\{ 0,2,4,\ldots,90\}$, calculated with 5 significant digits of precision. 
	In the case of $D=0$, the null block of the vector $\delta_{(0,0,0)}$ has been removed.}
 \label{table:Eigenvalues}
 \end{table}

Numerical approximation of the integrals $L_{{\bf k}}$ and $R_{{\bf k}}$ will rely on the following family of Bessel integrals for  $\sum_{j=1}^6 k_j=0$:
$$I_{{\bf k}}:= (2\pi)^{-5} \int_{(\mathbb{S}^1)^6} 	\ddirac{\sum_{j=1}^6 \omega_j}  \Big(\prod_{j=1}^6 \omega_j^{k_j} \d\sigma_{\omega_j}\Big)$$
$$= (2\pi)^{-1}\int_{\mathbb{R}^2} 	\prod_{j=1}^6 
J_{k_j}(|x|) \d x
= \int_0^\infty 	\prod_{j=1}^6 J_{k_j}(r) r\, \d r,$$
where the Bessel function $J_k$ is defined by 
$$ \int_{\mathbb{S}^1} \omega^k e^{-ix\cdot\omega} \d\sigma_\omega
= 2\pi (-i)^k J_k(|x|)(x/|x|)^k.$$
Indeed, writing
$$|\omega_4+\omega_5+\omega_6|^2-1=2+\sum_{j,k \in \{4,5,6\}, k\neq j} \omega_j\omega_k^{-1},$$
we obtain
$$L_{{\bf m},{\bf n}}=2I_{{\bf m},{\bf n}}+\sum_{\sigma\in S_3}
I_{{\bf m},{\bf n}+(1,-1,0)_\sigma},$$
$$R_{{\bf m},{\bf n}}=2I_{{\bf m}+{\bf n},(0,0,0)}+\sum_{\sigma\in S_3}
I_{{\bf m}+{\bf n},(1,-1,0)_\sigma}.$$ 

Using the fact that $J_{-k} = (-1)^{k} J_{k}$ and the above representation we see that $Q_{{\bf m}, {\bf n}} = Q_{-{\bf m}, -{\bf n}}$, so it suffices to consider $D \geq 0$.

To evaluate the integrals $I_{\bf k}$, we follow the 
scheme in  \cite{OST17}. 
We split the integrals into
\begin{equation}\label{eq:SplitI}
I_{\bf k}= 
 \int_0^S 	\prod_{j=1}^6 
J_{k_j}(r) r\, \d r+\int_S^R 	\prod_{j=1}^6 
J_{k_j}(r) r\, \d r+\int_R^\infty 	\prod_{j=1}^6 J_{k_j}(r) r\, \d r,
\end{equation}
with $S=3600$ and $R=63000$.
The first two integrals are evaluated with a Newton--Cotes 
quadrature rule. The step size is $0.003$ for the first integral and 
$0.05$ for the second integral.
In practice, this step involved tabulating the numerical values of $61$ Bessel functions at around $3\times 10^6$ points each,  with $20$ digit precision obtained via the software package {\it Mathematica} \cite{WM17}. This high precision lets the rounding errors be absorbed by the error estimates below.
 
The approximation error for the first integral in \eqref{eq:SplitI} was estimated  in \cite[\S 8]{OST17}, independently of the vector $\k$, by
\[
E_{\k,1} = 1.5\times 10^{-9}.
\]

The approximation error for the second integral in \eqref{eq:SplitI} was also estimated in \cite{OST17} by
\[
E_{\k,2} = C_{2} \prod_{j=1}^6 \Big(1+\frac {k_j^2} S\Big),
\]
where
\[
C_{2}
=
1.01^6(R-S)w^8\frac{6^3}5 \Big(\frac{2}{\pi(S-1)}\Big)^3\cosh^6(1) (R+1)
\]
with $S=3600$, $R=36000$ and $w=0.05$.
 
The third integral in \eqref{eq:SplitI} is approximated by analytic methods. 
Since $R=63000$ is large when compared to $n^2\leq 61^2$, we take advantage of the following well-known asymptotic formulae which are a simplified version of \cite[Corollaries 2.6 and 2.7]{OST17}.
\footnote{We record a typo on the first formula in \cite[Corollary 2.7]{OST17}, which should read as follows:
\begin{multline*}
\left|J_{0}^\pm(z)- \left(\frac{2}{\pi z}\right)^{\frac12}\cos(\omega_0)
-\frac 1{8z}  \left(\frac{2}{\pi z}\right)^{\frac12}\sin(\omega_0)\right|\le \\
\le\frac{9}{128|z|^2}\left(\frac{2}{\pi |z|}\right)^{\frac12} \cosh(|\Im (z)|)
\left(\frac{|z|}{|\Re(z)|}\right)^{\frac52}.
\end{multline*}}
\begin{lemma}\label{le:BesselAsymp}
	Let $\omega_n:=z-\frac{\pi}4-\frac{n\pi}2$ and $\hat{n}:=\max\{1,n\}$.
	If $n\geq 0$ and $z>\hat{n}^2$, then
	\begin{equation}\label{eq:BesselAsymp}
	\Big|J_n^{\pm}(z)-\Big(\frac2{\pi z}\Big)^{\frac12}\cos(\omega_n)\Big|
	\leq
	\Big(\frac2{\pi |z|}\Big)^{\frac12}\frac{\hat{n}^2}{|z|},
	\end{equation}
	\begin{multline}\label{eq:BesselAsymp2}
	\Big|J_n^{\pm}(z)-\Big(\frac2{\pi z}\Big)^{\frac12}\cos(\omega_n)+\frac{4n^2-1}{8z}\Big(\frac2{\pi z}\Big)^{\frac12}\sin(\omega_n)\Big|
	\leq
	\frac14\Big(\frac2{\pi |z|}\Big)^{\frac12}\frac{\hat{n}^4}{|z|^2}.
	\end{multline}\end{lemma}
\noindent Here the functions $J_n^\pm$ are obtained by writing $\cos(zt)=(e^{izt}+e^{-izt})/2$ in the Poisson integral representation for $J_n$,
and as such satisfy $J_n=(J_n^++J_n^-)/2$.
Using \eqref{eq:BesselAsymp}, we may split each Bessel function into a main term plus error.
Applying the distributive law yields one main integral of the form 
\begin{equation}\label{eq:ExactCosIntegral}
\int_R^\infty \Big(\frac2{\pi r}\Big)^3 \Big(\prod_{j=1}^6  \cos(\omega_{k_j})\Big) r\d r,
\end{equation}
which can be calculated exactly,
plus $2^6-1$ further terms involving one of the two factors $$\left( \frac{2}{\pi r} \right)^{\frac 1 2}\cos(\omega_{k_j}), J_{k_j}(r)-\left( \frac{2}{\pi r} \right)^{\frac 1 2}\cos(\omega_{k_j})$$
for each $j$. We call them cosine factor and error factor. 
For the main integral, observe that 
\[
\cos(r-\tfrac{\pi}4-\tfrac{k\pi}2)
=
(-1)^{\lfloor\frac{k}2\rfloor} \cdot \begin{cases}
\sin(r-\frac{\pi}4),& \text{if $k$ is odd},\\
\cos(r-\frac{\pi}4),& \text{if $k$ is even,}
\end{cases}
\]
and so \eqref{eq:ExactCosIntegral} equals a multiple of
$$\int_R^\infty \cos^6(r-\tfrac{\pi}4) r^{-2} \d r, \text{ or }\int_R^\infty \cos^4(r-\tfrac{\pi}4)\sin^2(r-\tfrac{\pi}4) r^{-2}\d r,$$
with sign determined by the parity of $\sum_{j=1}^3(\lfloor\frac{m_j}2\rfloor+ \lfloor\frac{n_j}2\rfloor)$.
{\it Mathematica} calculates these expressions with any prescribed accuracy.
For the further terms, consider first those consisting of an integral of a product of five cosine factors and one error factor. 

To estimate these six terms, we use the finer information given by \eqref{eq:BesselAsymp2}.
The sine term in \eqref{eq:BesselAsymp2} leads to integrals of the type
$$\frac{4m_1^2-1}{8} \int_R^\infty \Big(\frac{2}{\pi r}\Big)^3 \sin(\omega_{m_1})\cos(\omega_{m_2})\cos(\omega_{m_3})\cos(\omega_{\bf n})\d r$$
and similar terms with a different cosine factor replaced by a sine factor and corresponding prefactor.
The product of the six trigonometric functions is odd about the point $\frac{\pi}4$.
Thus the product integrates to 0 over each period.
On the period $[R+2\pi k, R+2\pi (k+1)]$, we may thus replace the weight $r^{-3}$ by the difference between itself and its mean over that interval, which in turn is bounded by $6\pi r^{-4}$.
Thus the sum of terms arising from the sine term in \eqref{eq:BesselAsymp2} is bounded by
$$E_{\k,3} = 3\pi \sum_{j=1}^6\hat{k}_j^2 \int_R^\infty \Big(\frac2{\pi}\Big)^3 r^{-4}\d r,$$
where $\hat{k}_j:=\max\{1,k_j\}$.
The sum of the six terms arising from the right-hand side of \eqref{eq:BesselAsymp2} can be estimated by
$$E_{\k,4} = \frac14\sum_{j=0}^6 \hat{k}_j^4 \int_R^\infty \Big(\frac2{\pi}\Big)^3 r^{-4}\d r.$$
Next come fifteen terms of the original $2^6-1$ terms which have four cosine factors and two error factors. It suffices to estimate these with
\eqref{eq:BesselAsymp}, since they benefit from an extra integration of a negative power of $r$. Their sum can be estimated by
$$E_{\k,5} = \sum_{i\neq j}\hat{k}_i^2 \hat{k}_j^2\int_R^\infty\Big(\frac2{\pi}\Big)^3 r^{-4}\d r,$$
where the sum is over ${6\choose 2}=15$ choices of distinct indices $i,j\in\{1,2,3,4,5,6\}$.

The remaining $2^6-1-6-15=42$ terms benefit from an integration of at least the negative fifth power of $r$, and are estimated even more crudely by
$$
E_{\k,6} = \sum_{i,j,\ell} \hat{k}_i^2 \hat{k}_j^2 \hat{k}_\ell^2 \int_R^\infty \Big(\frac2\pi\Big)^3 r^{-5}\d r
+\sum_{i,j,\ell,m} \hat{k}_i^2 \hat{k}_j^2 \hat{k}_\ell^2 \hat{k}_m^2 \int_R^\infty \Big(\frac2\pi\Big)^3 r^{-6}\d r$$
$$+\sum_{i,j,\ell,m,n} \hat{k}_i^2 \hat{k}_j^2 \hat{k}_\ell^2 \hat{k}_m^2 \hat{k}_n^2\int_R^\infty \Big(\frac2\pi\Big)^3 r^{-7}\d r
+\hat{k}_1^2 \hat{k}_2^2 \hat{k}_3^2 \hat{k}_4^2 \hat{k}_5^2 \hat{k}_6^2\int_R^\infty \Big(\frac2\pi\Big)^3 r^{-8}\d r,$$
where the sums are over tuples of distinct indices for a total of ${6\choose 3}=20$, ${6\choose 4}=15$, and ${6\choose 5}=6$ summands, respectively.

Addition of the error bounds $E_{\k,1}+\dotsb+E_{\k,6}$ yields error bounds for $I_{\k}$, which in turn give error bounds for the matrix coefficients $Q_{\mathbf{m},\mathbf{n}}$.
Applying Schur's test to each block with a fixed $D$ individually shows that the matrix consisting of the error bounds has operator norm less than $10^{-5}$.
These steps were again performed via {\it Mathematica}.
Since $10^{-5}$ is smaller than the minimal eigenvalues shown in Table~\ref{table:Eigenvalues}, we can conclude that the matrix $(Q_{\mathbf{m},\mathbf{n}})_{\mathbf{m},\mathbf{n}\in X}$ is positive definite.
This completes the proof of Theorem \ref{thm:main}.

\section{Further remarks}\label{sec:Remarks}

We conclude our discussion with several observations.

Table \ref{table:Eigenvalues} reveals that the smallest eigenvalues of the block $D$ is increasing in the parameter $D\ge 0$.
It might be interesting to find an analytic explanation of this fact.
  
\begin{figure}[htb]
\centering
\includegraphics[height=6cm]{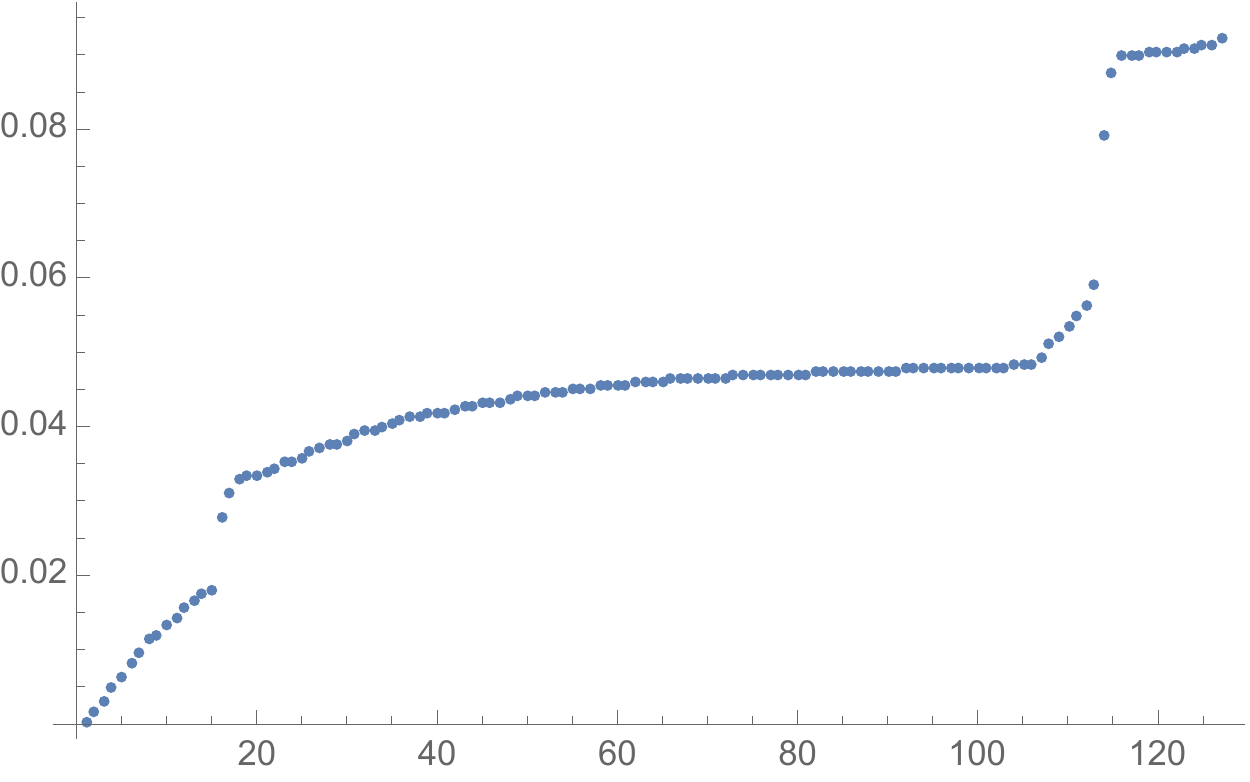}
\caption{ Plot of the eigenvalues 
  $0<\lambda_{1}\leq\lambda_{2}\leq\ldots\leq\lambda_{127}$ 
  of the approximation to the block $D=0$.}
\label{fig:EigenvaluesT}
\end{figure}

Zooming into the main block $D=0$, 
Figure \ref{fig:EigenvaluesT} shows the non-zero eigenvalues of this block. There is a cluster of very small eigenvalues. 
The corresponding eigenvectors seem to suggest that many of these small eigenvalues are related to functions on the circle that are mainly supported in neighborhoods of two antipodally symmetric points.
These functions are close competitors of constants for being maximizers.
A line of attack on this problem, say for larger  or infinite bandwidth, might be to understand how to analytically separate the effect of these antipodal functions.
The remaining eigenvalues may be sufficiently far from zero to allow for crude estimation.

We calculated the entries of the quadratic form $Q$
numerically. A look at these entries reveals some nice patterns such as circular structures, shown Figures \ref{fig:Circle30v1} and \ref{fig:Circle30v2} below. We do not know how to exactly describe or explain these structures independently of the numerical calculations. These structures merit further investigation.
 Each of the  six figures below shows a row of the block $D=0$.
This corresponds to fixing an index ${\bf m}_0$, and plotting the matrix entries corresponding to $Q_{{\bf m}_0,{\bf n}}$, where ${\bf n}$ ranges over all admissible values. Since $n_1+n_2+n_3=0$, we may
parametrize the entries of the row by $(n_1,n_2)$, which ranges
in a hexagonal region in $\Z^2$, shown in the figures as complement
of the shaded region.

\begin{figure}[htb]
\centering
\includegraphics[height=4.1cm]{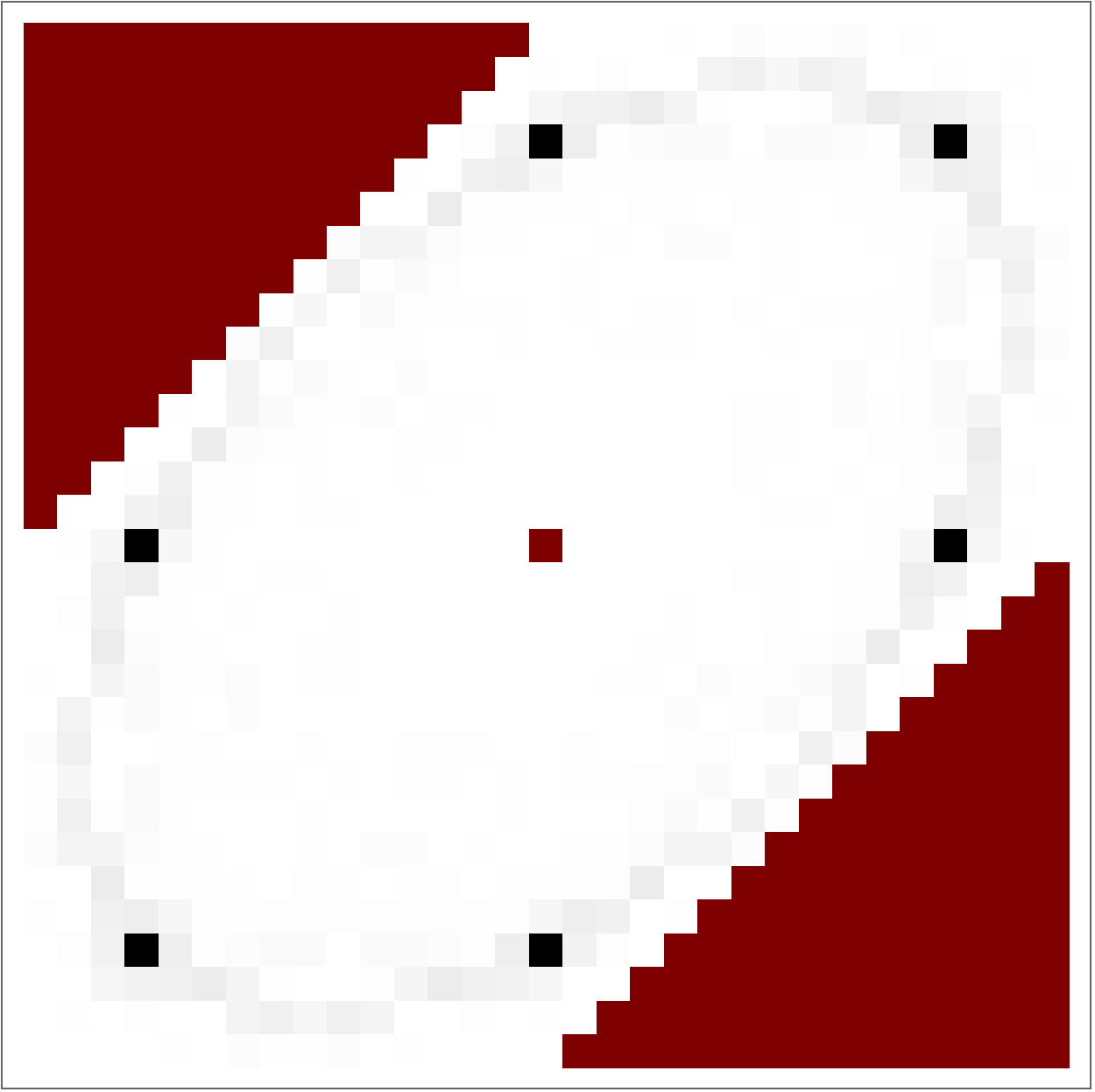}
\includegraphics[height=4.1cm]{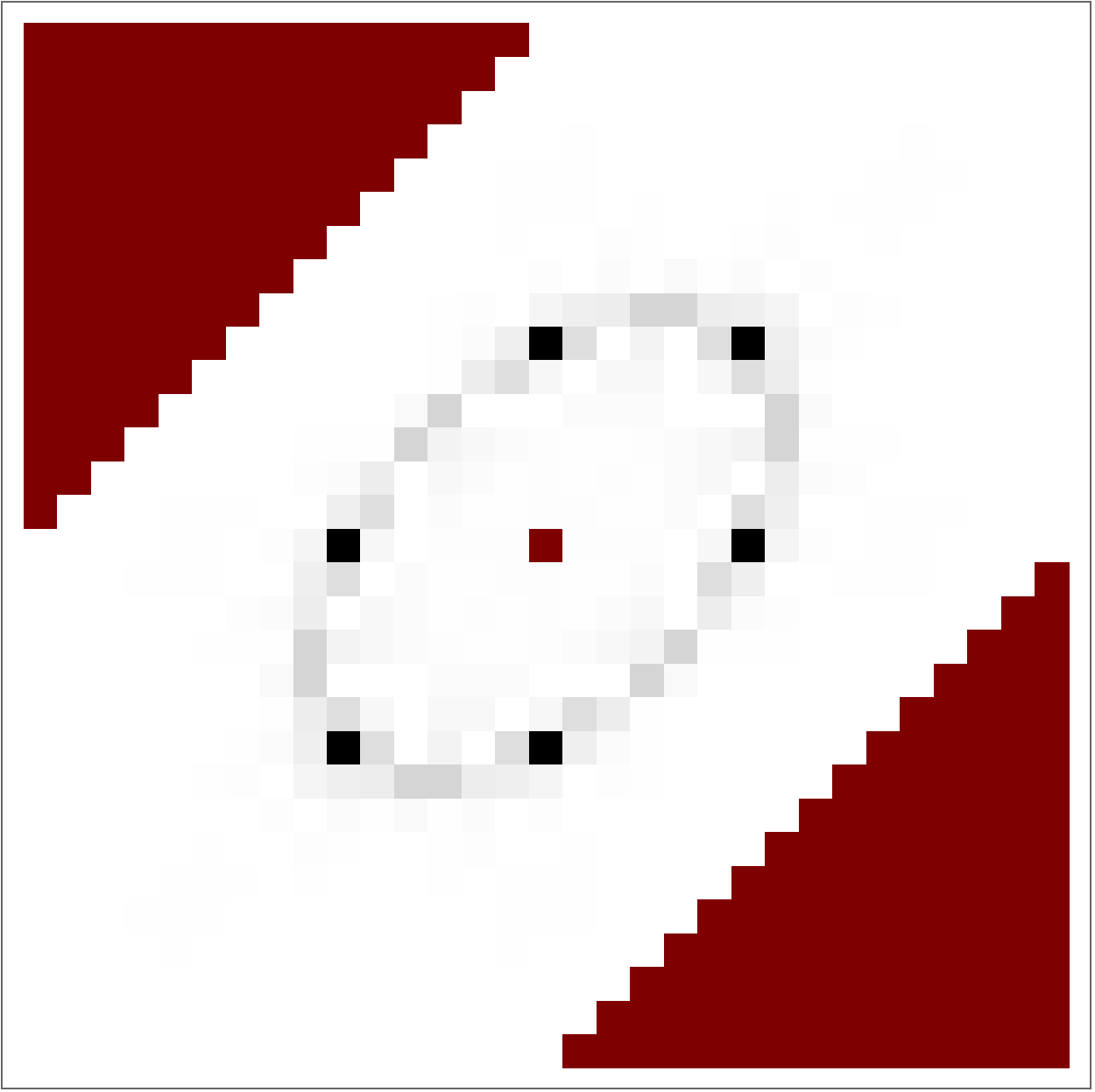}
\includegraphics[height=4.1cm]{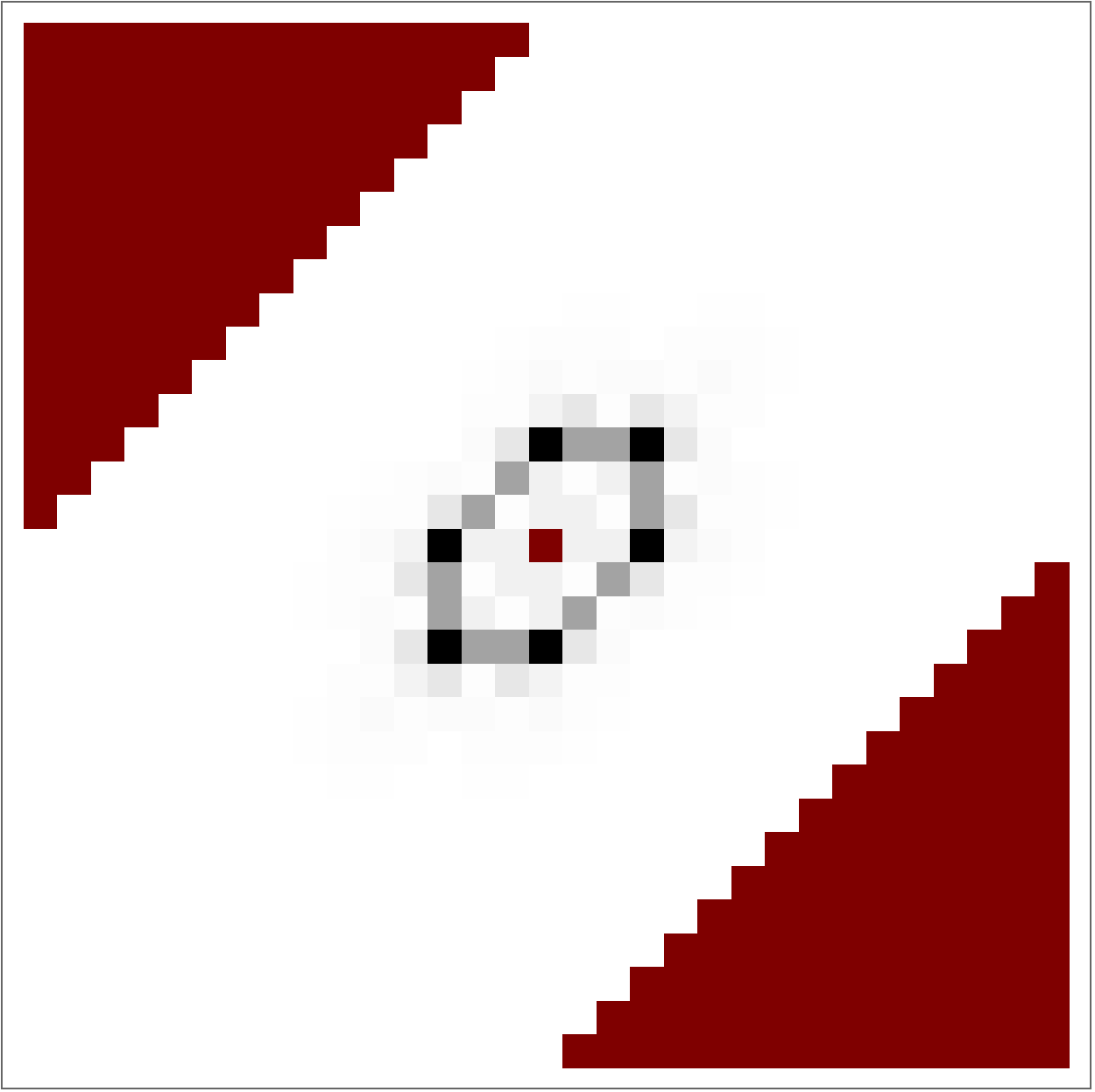}
\caption{${\bf m}_0=(-24,0,24)$, ${\bf m}_0=(-12,0,12)$,  ${\bf m}_0=(-6,0,6)$.}
\label{fig:Circle30v1}
\end{figure}

\begin{figure}[htb]
\centering
\includegraphics[height=4.1cm]{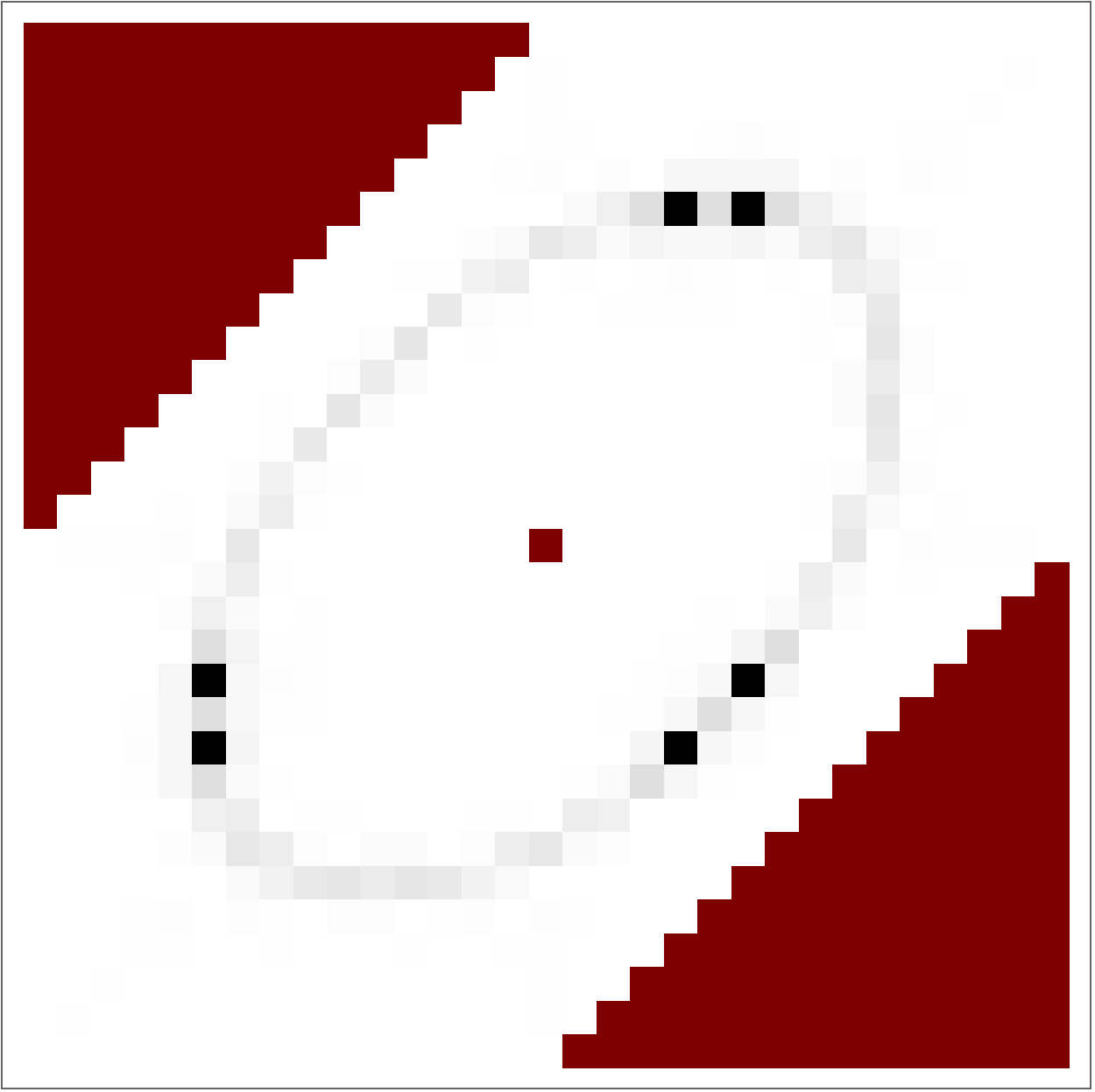}
\includegraphics[height=4.1cm]{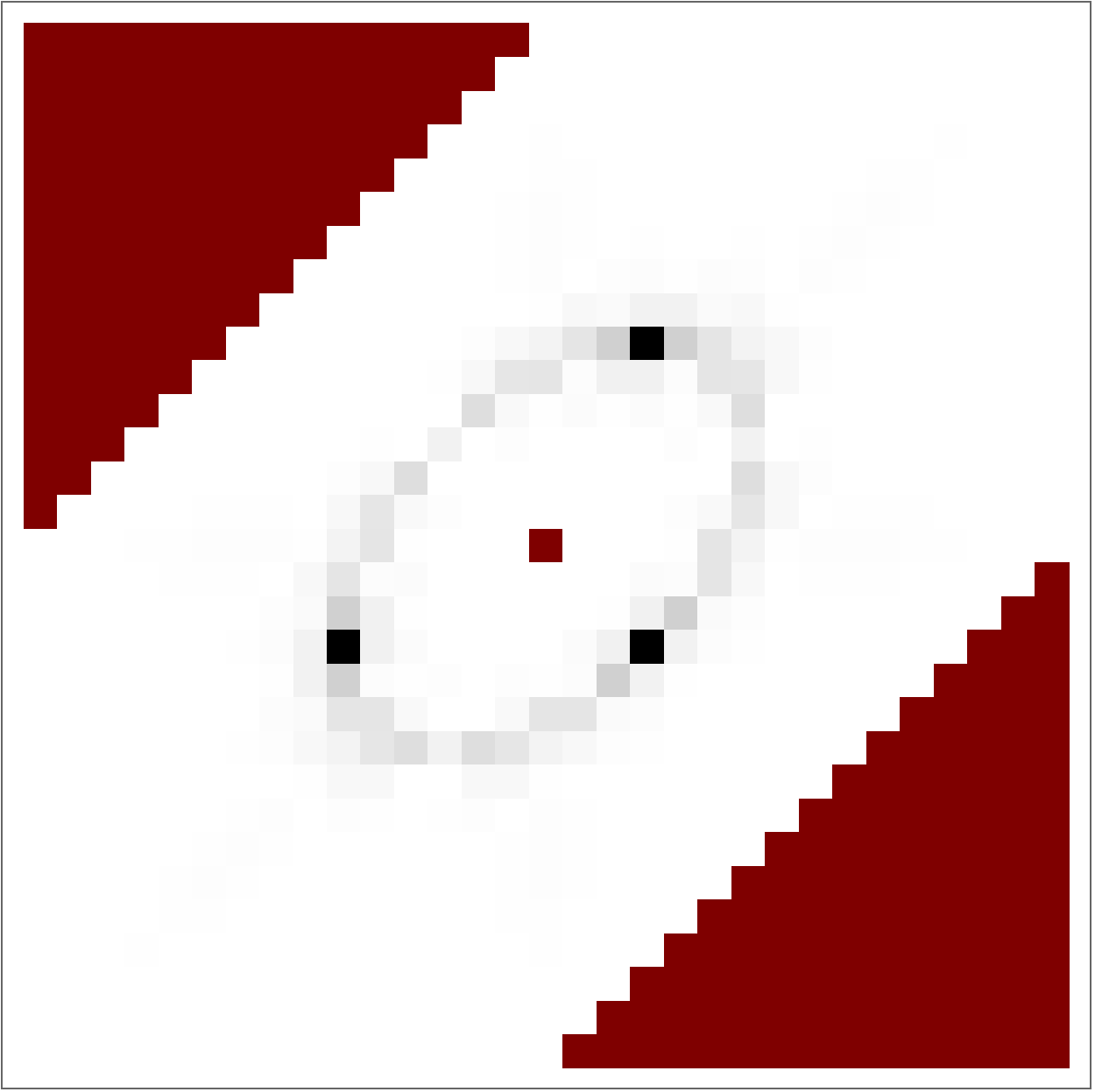}
\includegraphics[height=4.1cm]{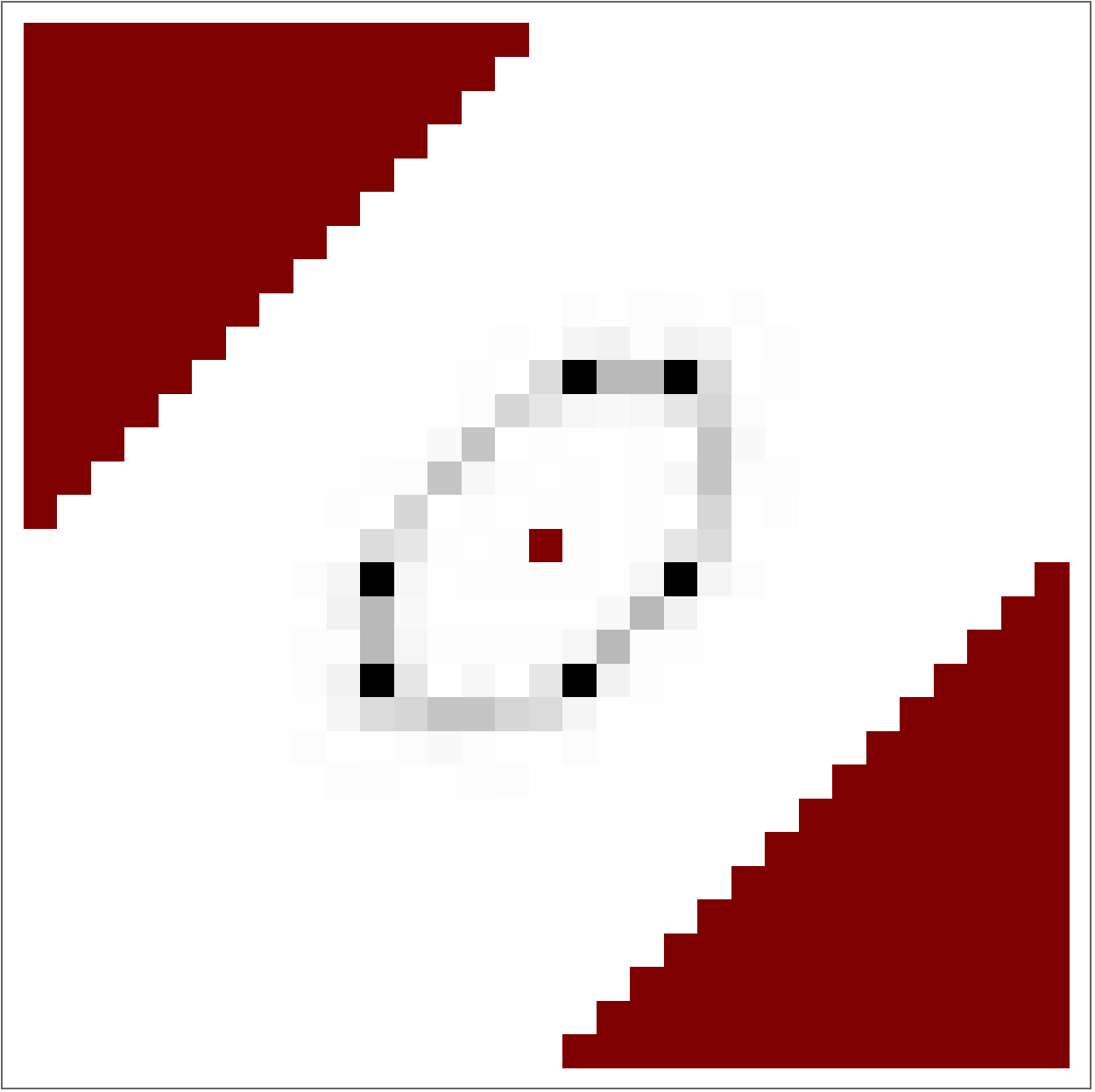}
\caption{${\bf m}_0=(-20,8,12)$, ${\bf m}_0=(-12,6,6)$, ${\bf m}_0=(-10,2,8)$.}
\label{fig:Circle30v2}
\end{figure}

We close with a remark on the central Bessel integral 
$$I_{(0,0,0,0,0,0)}=\int_0^\infty J_0^6(r)r\d r,$$
which up to a factor $(2\pi)^4$ is the conjectured sharp constant
$\Phi({\bf 1})$
in the Tomas--Stein adjoint restriction inequality. 
This integral appears in the following related context.
An {\it $n$-step uniform random walk} is a walk in the plane that starts at the origin and consists of $n$ steps of length 1 each taken into a uniformly random direction.

If $p_n$ denotes the radial density of the distance travelled after $n$ steps, then it is a classical result that $p_5(1)=I_{(0,0,0,0,0,0)}$, see e.g. \cite{BSWZ12}. In the same paper, the exact value of the integral
$$p_4(1)=\int_0^\infty J_0^5(r)r\d r=\frac{1}{2\pi^2}\sqrt{\frac{\Gamma(\frac1{15})\Gamma(\frac2{15})\Gamma(\frac4{15})\Gamma(\frac8{15})}{5\Gamma(\frac7{15})\Gamma(\frac{11}{15})\Gamma(\frac{13}{15})\Gamma(\frac{14}{15})}}$$
is determined resorting to  striking modularity properties of the function $p_4$, see \cite[Theorems 4.9 and 5.1]{BSWZ12}.
The corresponding problem with a sixth power remains open.

\section*{Acknowledgements}
We thank Emanuel Carneiro, Damiano Foschi and Felipe Gon\c{c}alves for stimulating discussions. 
The software {\it Mathematica} was used to perform the numerical tasks described in \S\ref{sec:Numerics} and \S\ref{sec:Remarks}. 
The authors acknowledge support by the Hausdorff Center for Mathematics and the Deutsche Forschungsgemeinschaft through the Collaborative Research Center 1060.

\end{document}